\documentclass[preprint,12pt]{elsarticle}

\usepackage{graphicx}
\usepackage{comment}
\usepackage{amssymb}
\usepackage{lineno}
\usepackage{hyperref}
\usepackage{amsmath,amsfonts,amsthm,amssymb,graphicx,subfigure,mathtools,tikz,hyperref,bm,blindtext,float}
\usepackage{diagbox}
\usepackage{geometry}
\usepackage[ruled]{algorithm2e} 

\journal{Journal Name}

\begin{document}

\begin{frontmatter}

%% Title, authors and addresses

%\title{Multi-fidelity physics-informed neural networks for function approximations and inverse PDE problems}
\title{Physics-informed neural networks with residual/gradient-based adaptive sampling methods for solving PDEs with sharp solutions}

%% use the tnoteref command within \title for footnotes;
%% use the tnotetext command for the associated footnote;
%% use the fnref command within \author or \address for footnotes;
%% use the fntext command for the associated footnote;
%% use the corref command within \author for corresponding author footnotes;
%% use the cortext command for the associated footnote;
%% use the ead command for the email address,
%% and the form \ead[url] for the home page:
%%
%% \title{Title\tnoteref{label1}}
%% \tnotetext[label1]{}
%% \author{Name\corref{cor1}\fnref{label2}}
%% \ead{email address}
%% \ead[url]{home page}
%% \fntext[label2]{}
%% \cortext[cor1]{}
%% \address{Address\fnref{label3}}
%% \fntext[label3]{}

%% use optional labels to link authors explicitly to addresses:
%% \author[label1,label2]{<author name>}
%% \address[label1]{<address>}
%% \address[label2]{<address>}

\author[Xiamen]{Zhiping Mao}
\ead{zpmao@xmu.edu.cn}
\author[hust]{Xuhui Meng\corref{cor}}
% \author[brown,PNNL]{George Em Karniadakis\corref{cor}}
%\cortext[cor1]{Corresponding author}
%\author[brown]{George Em Karniadakis\corref{cor}}
\cortext[cor]{Corresponding author}
\ead{xuhui\_meng@hust.edu.cn}

\address[Xiamen]{School of Mathematical Sciences, Fujian Provincial Key Laboratory of Mathematical Modeling and High-Performance Scientific Computing, Xiamen University, Xiamen, China}
\address[hust]{Institute of Interdisciplinary Research for Mathematics and Applied Science, School of Mathematics and Statistics, Huazhong University of Science and Technology, Wuhan, 430074, China}
% \address[PNNL]{Pacific Northwest National Laboratory, Richland, WA, 99354 USA}

\begin{abstract}
%% Text of abstract
We consider solving the forward and inverse PDEs which have sharp solutions using physics-informed neural networks (PINNs) in this work. In particular, to better capture the sharpness of the solution, we propose adaptive sampling methods (ASMs) based on the residual and the gradient of the solution. 
We first present a residual only based ASM algorithm denoted by ASM I. In this approach, we first train the neural network by using a small number of residual points and divide the computational domain into a certain number of sub-domains, we then add new residual points in the sub-domain which has the largest mean absolute value of the residual, and those points which have largest absolute values of the residual in this sub-domain will be added as new residual points. We further develop a second type of ASM algorithm (denoted by ASM II) based on both the residual and the gradient of the solution due to the fact that only the residual may be not able to efficiently capture the sharpness of the solution. The procedure of ASM II is almost the same as that of ASM I except that in ASM II, we add new residual points which not only have large residual but also large gradient. To demonstrate the effectiveness of the present methods, we employ both ASM I and ASM II to solve a number of PDEs, including Burger equation, compressible Euler equation, Poisson equation over an L-shape domain as well as high-dimensional Poisson equation. It has been shown from the numerical results that the sharp solutions can be well approximated by using either ASM I or ASM II algorithm, and both methods deliver much more accurate solution than original PINNs with the same number of residual points. Moreover, the ASM II algorithm has better performance in terms of accuracy, efficiency and stability compared with the ASM I algorithm. This means that the gradient of the solution improves the stability and efficiency of the adaptive sampling procedure as well as the accuracy of the solution. 
Furthermore, we also employ the similar adaptive sampling technique for the data points of boundary conditions if the sharpness of the solution is near the boundary. The result of the L-shape Poisson problem indicates that the present method can significantly improve the efficiency, stability and accuracy. 
\end{abstract}

\begin{keyword}
physics-informed neural networks \sep adaptive sampling \sep high-dimension \sep L-shape Poisson equation \sep accuracy
%% keywords here, in the form: keyword \sep keyword

%% MSC codes here, in the form: \MSC code \sep code
%% or \MSC[2008] code \sep code (2000 is the default)

\end{keyword}

\end{frontmatter}

%%
%% Start line numbering here if you want
%%
%\linenumbers

%% main text
\section{Introduction}
\label{intro}
Deep learning methods (e.g., deep neural network (DNN), etc) have recently achieved great progress in an extensive fields including speech recognition, natural language processing, image recognition analysis, and so on \cite{lecun2015deep,mikolov2011strategies,hinton2012deep,sainath2013deep,krizhevsky2012imagenet,tompson2014joint}. Most recently, the DNN has also been extended to encode the physical laws (e.g., mass/momentum/energy conservation law) in the form of partial differential equations (PDEs), which yields the physics-informed neural networks (PINNs) \cite{raissi2019physics}, to solve PDEs. See also \cite{yu2018deep,han2018solving,long2018pde,long2019pde,sirignano2018dgm,pang2019fpinns,meng2020composite,raissi2019physics,zhang2019quantifying,lu2021deepxde} and references therein. The PINN can solve both the forward problems of PDEs and discover the unknown parameters (mathematical modelling). And it has become a powerful tool in various applications, such as fluid mechanics \cite{raissi2018hidden,raissi2019deep, cai2021physics}, hydrogeophysics \cite{meng2020composite}, transport phenomena in porous media \cite{pang2019fpinns}, and compressible Euler equations~\cite{mao2020physics, JAGTAP2022111402}, since the PINN algorithm has several advantages compared with the classical numerical methods: (1) the algorithm of PINN is easy coding with the help of auto-differentiation~\cite{A.Baydin2018JMLR}; (2) it is a meshless method which can easy deal with complex geometric problems; (3) PINNs can also solve high-dimensional problems without suffering the issue of curse of dimensionality; (4) moreover, it is efficient to solve the inverse problems (e.g., learning the model parameters). 

An important step and a frequency asked question for PINNs is about the sampling of the residual points. The distribution of the residual points has a significant influence to the efficiency and accuracy of the PINN algorithm. For instance, it has been shown in \cite{mao2020physics} that a cluster distribution using more residual points around the discontinuity results in more accurate solution. However the cluster distributions of the residual points used in \cite{mao2020physics} are based on a prior knowledge of the solution which is usually not available. To resolve this issue, several affects have been made by using adaptive methods similar as the adaptive mesh refinement methods used in the conventional numerical methods for problems with discontinuous solutions \cite{berger1984adaptive, berger1989local, verfurth1994posteriori, baeza2006adaptive}. A basic idea of the adaptive mesh refinement methods is to refine the grids near the discontinuities based on the gradient of the numerical solutions, which can significantly enhance the accuracy as well as reduce the computational cost. 
To this end, a residual-based adaptive refinement method is proposed in \cite{lu2021deepxde} and later a similar refinement method~\cite{yu2022gradient} is proposed based on both the residual and the gradient of the residual. Recently, a kind of sampling method using the probability density function related to the residual is developed in \cite{wu2023comprehensive}. Similarly, from the point view of probability, Zhou et. al. developed an failure-informed self-adaptive sampling method using failure probability based indicator in \cite{gao2022failure}. See also \cite{guo2022novel,han2022residual} adaptive causal sampling method and weighted adaptive sampling method.

The aforementioned adaptive sampling methods are all residual based, namely, these adaptive methods only use the information of the residual. However, the gradient of the (predicated) solution is also an important and potential estimator in detecting the discontinuity and then determining the new residual points. Motivated by this, we proposed in this paper adaptive sampling methods (ASMs) based on the residual and/or the gradient of the solution for the problems which have sharp solutions. More precisely, we develop a residual only based ASM and an ASM based on both the residual and the gradient of the solution. We employ both ASMs to several linear and nonlinear equations including Burgers equation, compressible Euler Equation, Poisson equation over an L-shape domain as well as high-dimensional Poisson equation. All numerical results indicate that the residual/gradient-based method has better performance than the residual-only based method in terms of stability, efficiency and accuracy. Moreover, for the problem whose solution has sharpness near the boundary, we also employ similar adaptive technique for the boundary data points resulting in more accurate and stable solutions.

The rest of the paper is organized as follows: The adaptive sampling algorithms are presented in Sec. \ref{meth}, with several examples in Sec. \ref{results} to demonstrate its effectiveness. Finally, a summary with some discussions are displayed in Sec. \ref{summary}.

\section{Methodology}
\label{meth}
\begin{figure}
    \centering
    \includegraphics[width = 0.7\textwidth]{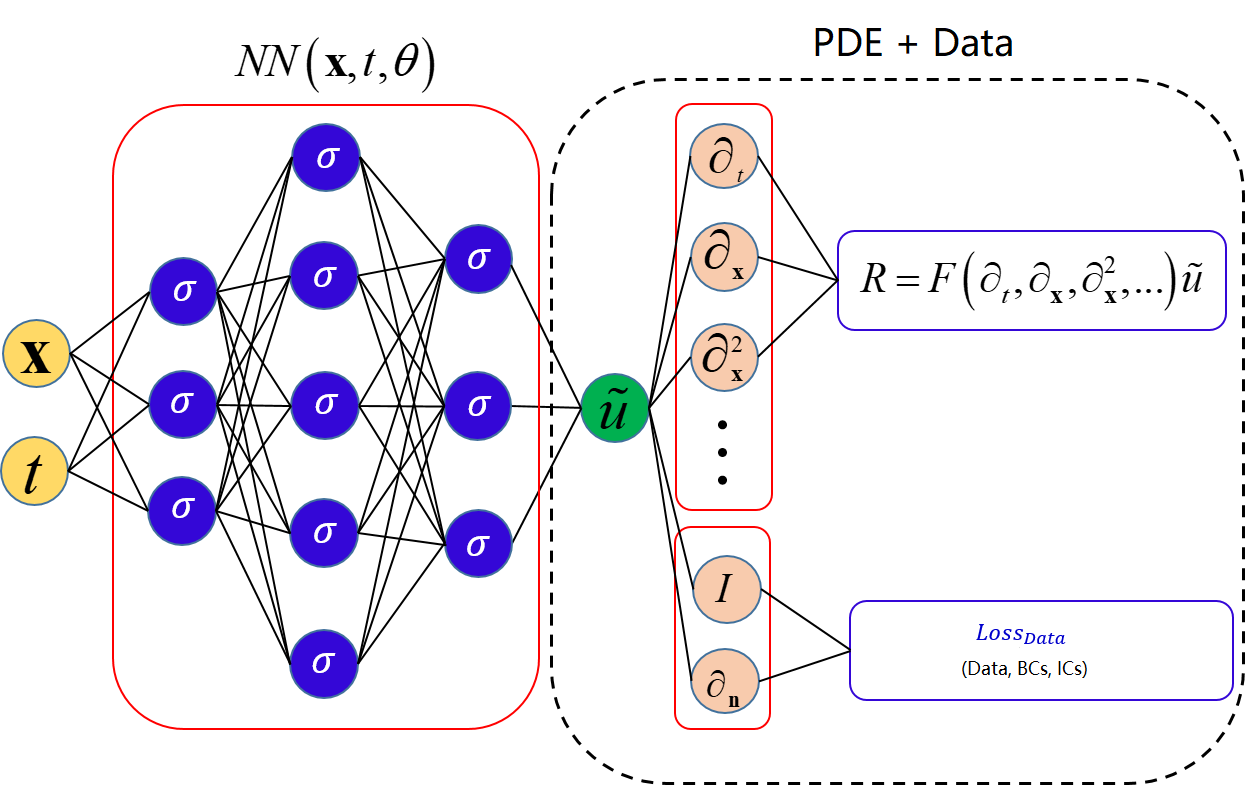}
    \caption{Schematic of PINNs. Yellow dots: inputs of the NN, blue dots: activation function, green dots $\tilde{u}$: output of the neural network (which is used to approximate the solution of the PDE), orange dots: differential operator, $R$ represents the residual of the PDE, $Loss_{Data}$ represents the mismatch of data, the data can be boundary/initial conditions or any other available data. 
    %$u$ is the solution for a PDE, which is approximated by the first NN i.e., $u \approx $ NN$(\bm{x}, t, \theta)$ = $\tilde{u}$, $\theta$ denotes the hyperparameter in the first NN.
    }\label{pinns}
\end{figure}

We present in this section the basic idea of ASM and consequently two ASM algorithms. We begin by introducing the algorithm of PINNs. 
In PINNs, the loss consists of two parts: the loss of the residual and the loss of data, namely,
$$\mbox{Loss} = \omega_R \mbox{Loss}_{R} + \omega_{D}\mbox{Loss}_{data},$$
where $\mbox{Loss}_{R}$ is the mean square of the residuals for PDEs defined by 
$$\mbox{Loss}_{R} = \frac{1}{N_R}\sum^{N_R}_i R_i^2,$$
and $\mbox{Loss}_{data}$ is corresponding to the data defined by 
$$\mbox{Loss}_{data} = \frac{1}{N_{D}}\sum^{N_{D}}_i (u^{data}_i - u^{NN}_i)^2,$$
the data can be the initial/boundary conditions or any other available data,
here $N_{\psi}$ and $\omega_{\psi}$ ($\psi = R, \; D$) are the number and weight for each component, respectively.
A schematic of PINNs is shown in Fig. \ref{pinns}.

We then introduce how to adaptively add residual points based on the residual and the gradient of the solution. We first introduce the residual only-based adaptive sampling method denoted by ASM I, The ASM I is similar as the residual-based adaptive refinement method proposed in \cite{lu2021deepxde} except that in ASM I, we use domain decomposition and compute the mean value of the residual for each sub-domain, i.e., $\mathcal{E}_i, i=1,\ldots, {\bm S}$, and then obtain the new residual points in the sub-domain which has the largest absolute value of the residuals. The detail of the algorithm is given in Algorithm \ref{alg:ASMI}. 
We point out here that the use of the mean value of the residual for each sub-domain can eliminate the occasionality of bad residual points and make the iteration procedure of the adaptive sampling step to be more stable. 
%Moreover, it can improve the performance of the present 

\begin{algorithm}[htbp]
%    \SetAlgoNoLine  %去掉之前的竖线 
\caption{ASM I: Residual only-based adaptive sampling method}\label{alg:ASMI}
\begin{enumerate}
\item[Step 1] Train the PINNs based on the residual point set $\mathcal{T}$;
\item[Step 2] Divide the computational domain into $\bm{S}$  subdomains $\Omega_i, i =1,\ldots, {\bm S}$, and then compute the mean value  ($\mathcal{E}_i, i=1,\ldots, {\bm S}$) of the absolute value of the residual ($R$), namely, 
$$\mathcal{E}_i = \int_{\Omega_i} |R| d{\bm x},$$
using the Gaussian quadrature method for each sub-domain (See Fig. \ref{integral} for a 1D representative case).
\item[Step 3]  Let $\mathcal{E}_{max}$ denote the largest mean value of the PDE residuals for all subdomains, i.e., $\mathcal{E}_{max} = \max_{i = 1, \ldots {\bm S}}\{\mathcal{E}_i\}$. If $\mathcal{E}_{max} < \mathcal{E}_c$, iteration stops ($\mathcal{E}_c$ is a user-defined convergence criterion). Otherwise, randomly sample $\mathcal{M}$ data in the sub-domain  which has the largest $\mathcal{E}$, and then compute the residuals for the selected points. Add $m_{R}$ points which have the largest absolute value of the residuals in the residual point set $\mathcal{T}$.
%\item (4) For discontinuous problems, e.g., shock wave, to resolve the discontinuity accurately and efficiently, we can also add $\mathcal{P}$ points where $\nabla \tilde{u}$ are the largest after $\mathcal{E} < \mathcal{E}_g$. 
\item[Step 4] Repeat Step 1 - 3.
\end{enumerate}
\end{algorithm}

\begin{figure}
    \centering
    \includegraphics[width = 0.6\textwidth]{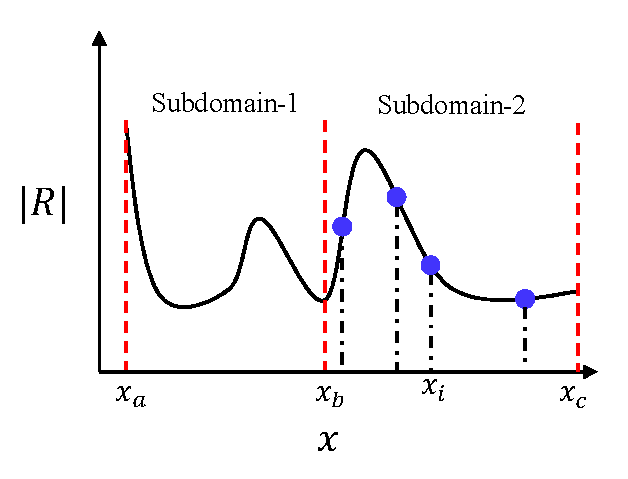}
    \caption{Schematic for the computation of mean value of the absolute value ($\mathcal{E}$) in a sub-domain. Based on the Gaussian quadrature, the $\mathcal{E}$ is calculated as $\mathcal{E} = \frac{1}{x_c - x_b} \int^{x_c}_{x_b} |R(x)| dx \approx \frac{1}{x_c - x_b} \sum^N_i \omega_i |R(x_i)|$, where $x_i$ is the Gauss quadrature point (blue dot),  $\omega_i, i=1,\ldots, N,$ are the Gauss quadrature weights, $N$ is the number of the quadrature points.}
    \label{integral}
\end{figure}

As mentioned previously, Only the residual may be not good enough to find good new added residual points, especially for the problems that have steep solutions. As we know, the gradient of the solution can be used to measure the sharpness of the solution. Therefore, we proposed a new adaptive sampling based on both the residual and the gradient of the solution. The procedure is similar as the that proposed in Algorithm \ref{alg:ASMI} except that we add new $m_{R}$ and $m_{\nabla u}$ residual points where the residual and the gradients of the solution are large, respectively. See the detail in Algorithm \ref{alg:ASMII}. If $m_{\nabla u} = 0$, it degenerates to the ASM I algorithm.

\begin{algorithm}[htbp]
%    \SetAlgoNoLine  %去掉之前的竖线 
\caption{ASM II: Residual/gradient-based adaptive sampling method}\label{alg:ASMII} 
\begin{enumerate}
\item[Step 1] Train the PINNs based on the residual point $\mathcal{T}$;
\item[Step 2] Divide the computational domain into $\bm{S}$  subdomains $\Omega_i, i =1,\ldots, {\bm S}$, and then compute the mean value  ($\mathcal{E}_i, i=1,\ldots, {\bm S}$) of the absolute value of the residual ($R$) and the mean value  ($\mathbf{E}_i, i=1,\ldots, {\bm S}$) of the absolute value of the gradient of the solution ($\nabla u$), namely, 
$$\mathcal{E}_i = \int_{\Omega_i} |R| d{\bm x}, \quad
\mathbf{E}_i = \int_{\Omega_i} |\nabla u| d{\bm x}
$$
using the Gaussian quadrature method for each sub-domain; 
\item[Step 3] Let $\mathcal{E}_{max} = \max_{i = 1, \ldots {\bm S}}\{\mathcal{E}_i\}$ and $\mathbf{E}_{max} = \max_{i = 1, \ldots {\bm S}}\{\mathbf{E}_i\}$. If $\mathcal{E}_{max} < \mathcal{E}_c$, iteration stops. Otherwise, randomly sample $\mathcal{M}$ data in the sub-domain which has the largest $\mathcal{E}$ and $\mathcal{M}$ data in the sub-domain which has the largest $\mathbf{E}$, and then compute the residuals and gradients for selected points, respectively. Add $m_{R}$ points which have the largest $|R|$ and $m_{\nabla u}$ points which have the largest $|\nabla u|$ in the residual point set $\mathcal{T}$.
%\item (4) For discontinuous problems, e.g., shock wave, to resolve the discontinuity accurately and efficiently, we can also add $\mathcal{P}$ points where $\nabla \tilde{u}$ are the largest after $\mathcal{E} < \mathcal{E}_g$. 
\item[Step 4] Repeat Step 1 - 3.
\end{enumerate}
\end{algorithm}

\section{Results and discussions}\label{results}
To demonstrate the effectiveness of the present methods, we present in this section several numerical examples for the problems which have sharp solutions. In particular, we consider the forward and inverse problems of Burgers equation in subsection \ref{sec:burgers:for} and subsection \ref{sec:burgers:inv}, respectively. We then solve the compressible Euler equation in subsection \ref{sec:euler}, and the Poisson equation over an L-shape domain and the high-dimensional Poisson equation in subsection \ref{sec:Poisson:Lshape} and subsection \ref{sec:Poisson:high}, respectively.

\subsection{Burgers equation}\label{sec:burgers:for}
We begin by considering the forward problem of Burgers equation, namely, we shall solve the following equation 
\begin{equation}
\begin{split}
    &\partial_t u + u\partial_x u = \nu \partial^2_x u, ~x \in(-1, 1), ~ t \in [0, 1],\\
    &u(x, 0) = -\sin(\pi x), \\
    &u(-1, t) = u(1, t) = 0
\end{split}
\end{equation}
by using the present methods, 
here $u$ is the fluid velocity, and $\nu = 0.01/ \pi$ corresponds to the kinematic viscosity. 
It is well know that though the initial condition is smooth, the solution of Burgers equation becomes sharp at $x= 0$ as time evolves.

We first employ 2000 randomly sampled residual points and 100 randomly sampled initial points and 200 randomly sampled boundary points for the initial training, and the training points will be added according to the present ASM algorithms.
The architecture of the NNs is $4\times 40$ ($w \times l$, $w$ and $l$ denote the depth and width of the DNN, respectively). At each iteration, we first train the neural network using Adam optimizer and then transfer to the L-BFGS optimizer if the loss less than $2\times 10^{-3}$ or the number of epoches of Adam optimizer is achieve 10000. The learning rate for the Adam is $8\times 10^{-4}$. The computational domain is divided into $x \times t = 40 \times 20$ uniform subdomains, and the number of the Gaussian quadrature points in each domain is set as $4 \times 4$. In addition, $\mathcal{M} = 80$, $\mathcal{E}_c = 0.5\%$, and we set $m_{R} = 2$ for ASM I while set $m_{R} = 1, \; m_{\nabla u} =1$ for ASM II. 

\begin{figure}
    \centering
    \subfigure[]{\label{burgersa}
    \includegraphics[width=0.8\textwidth]{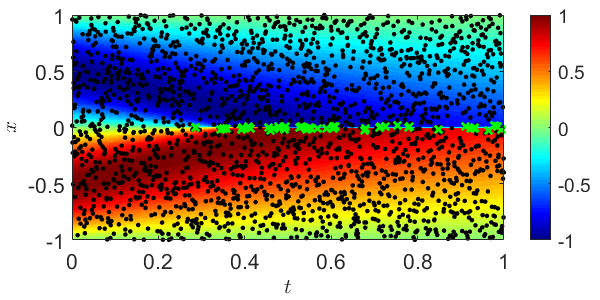}}
    \subfigure[]{\label{burgersb}
    \includegraphics[width=0.45\textwidth]{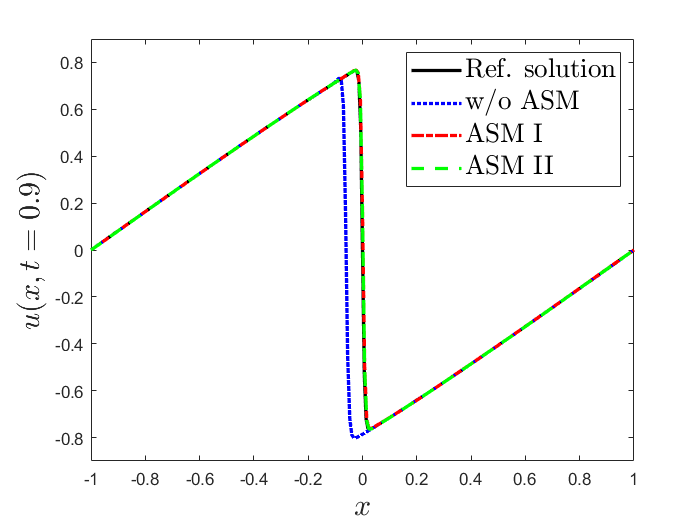}}
    \subfigure[]{\label{burgersc}
    \includegraphics[width=0.45\textwidth]{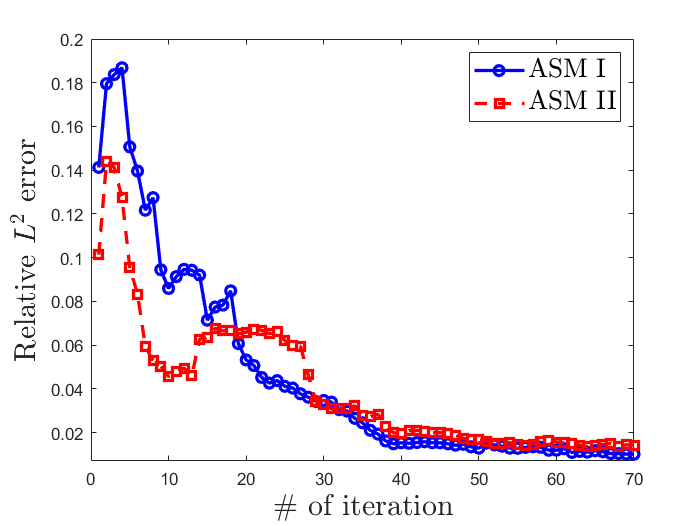}}
    \caption{\label{burgers}
    Results for the forward problem of Burgers equation. For ASM I, $m_{R} = 2$ while for ASM II, $m_{R} = 1, \; m_{\nabla u} =1$.
    (a) A representative distribution of the residual points  for solving the Burgers equation. Black dots: initial residual points, Green cross: added residual points. 
    (b) Comparison between the solutions from the present method and the reference solution. Dotted line: PINN solution using 2140 random residual points without using the ASM, Dash-dotted line: PINN solution using  ASM I.
    (c) Convergence the relative $L^2$ error using ASM I (blue solid line with circle) or ASM II (red dash line with square). The results are the mean values of the relative $L^2$ errors obtained from 10 runs.}
\end{figure}

We found that the distributions of the added points by using ASM I and ASM II are quite similar, and then we only show the distribution of the added points for the case of ASM I in Fig. \ref{burgersa}. It is interesting to find that almost of all the added residual points locate around $x = 0$ (green cross in Fig. \ref{burgersa}), and it clearly enhance the predictive accuracy for the PINNs. The predictions of $u$ at $t = 0.9$ are shown in  Fig. \ref{burgersb} for both the cases of ASM I and ASM II. For the purpose of comparison, we also randomly select 2140 residual points and train the PINNs without using the ARM, see the blue dotted line in Fig. \ref{burgersb}. Clearly, we observe that the solutions obtained by using ASM I and ASM II agree quite well with the reference solution while there is a big gap between the solution obtained without ASM and the reference solution.
In addition, we show the convergence of the relative $L^2$ error for $u$ showing that the relative $L^2$ error decreases asymptotically as the number of added residual points increases.
%We  note that the total number of the residual points is 2044. For the purpose of comparison, we randomly select 2044 residual points and train the PINNs without using the ARM. We observe from Figs. \ref{burgersb}-\ref{burgersc} that the PINNs with the ARM deliver more accurate predictions.

We further proceed to investigate the effect of the number of subdomains and the quadrature points on the predictive accuracy. We found that larger numbers of the subdomains and quadrature points given better prediction. We present one case with ASM I (here we set $m_{R}=1$) to demonstrate the effect of the number of subdomains and the quadrature points. We test a set of randomly sampled initial residual points with two different setups for the subdomains and quadrature points: (1) $x \times t = 20 \times 10$ uniform sub-domains with $x \times t = 4 \times 4$ quadrature points in each sub-domain, (2) $x \times t = 20 \times 10$ uniform sub-domains with $x \times t = 5 \times 5$ quadrature points in each sub-domain. The result is shown in Table \ref{burgers_subpoints} showing that increasing the number of the subdomains and quadrature points would improve the predictive accuracy.
%, which is the same as the results obtained in Sec. \ref{sec_poisson}.

\begin{table}[htbp]
\centering
 \begin{tabular}{c|ccc}
  \hline \hline
  ~ & $\{(20 \times 10), (4 \times 4)\}$ & $\{(20 \times 10), (5 \times 5)\}$ & $\{(40 \times 20), (4 \times 4)\}$  \\ \hline
 Added points    &  14 & 64 & 44   \\
 Relative $L^2$ error    &  $2.02 \%$ & $1.29 \%$ & $0.68 \%$  \\
  \hline \hline
 \end{tabular}
  \caption{\label{burgers_subpoints}
  Results for the forward problem of the Burgers equation using ASM I ($m_{R} = 1$): Relative $L^2$ errors for $u$ using different numbers of subdomains and quadrature points. The numbers in the first parenthesis of the brace represent the number of subdomains in the form of $x \times t$, and those of the second parenthesis are the numbers of quadrature points in each subdomain ($x \times t$).}
\end{table}

\subsection{Inverse problem of Burgers equation}\label{sec:burgers:inv}
In the previous example, we found that the results obtained by using ASM I and ASM II are quite similar to each other from the point view of stability, accuracy and efficiency. To investigate how the gradient of the solution affects the adaptive sampling, we now consider the inverse problem of Burgers equation, i.e., we would like to learn the value of the viscosity $\nu$ with data of $u$, by comparing the ASM II with ASM I.

Initially, we use 1000 randomly distributed residual points and 1000 randomly distributed data as well as 100 randomly sampled points for IC and 200 randomly sampled points for both BCs.
The architecture of the NNs is the same as the forward problem of Burgers equation, i.e., $4\times 40$. At each iteration, we first train the neural network using Adam optimizer with max epoch 8000. The Adam optimizer stops if the loss is less than $10^{-3}$. If the number of epochs using Adam is less than 1000, the L-BFGS optimizer is then activated. The learning rate for the Adam is $10^{-3}$. The computational domain is divided into $x \times t = 40 \times 20$ uniform subdomains, and the number of the Gaussian quadrature points in each domain is set as $4 \times 4$, $\mathcal{M} = 80$, $\mathcal{E}_c = 0.9\%$, and we set $m_{R} = 2$ for ASM I while set $m_{R} = 1, \; m_{\nabla u} =1$ for ASM II. The initial value of $\nu$ is set to be $0.02/\pi$.

\begin{figure}[htbp]
    \centering
    \includegraphics[width = 0.6\textwidth]{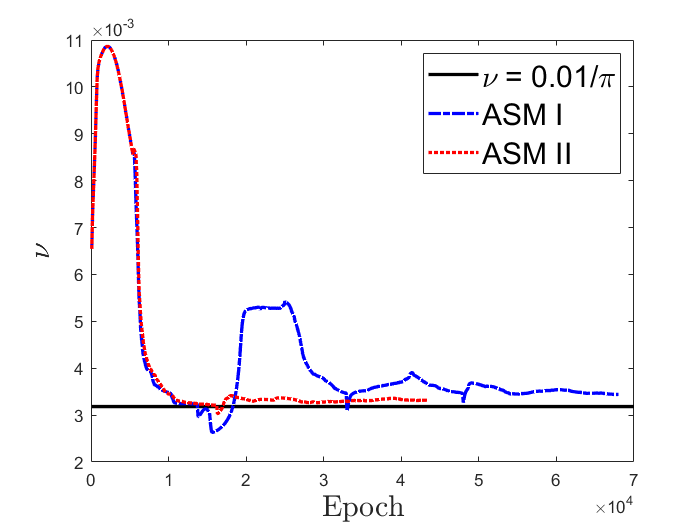}
    \caption{Results for the inverse problem of Burgers equation. For ASM I, $m_{R} = 2$ while for ASM II, $m_{R} = 1, \; m_{\nabla u} =1$: A representative comparison of the convergence of the learned $\nu$ between the cases of ASM I and ASM II, showing that the learned value of $\nu$ obtained by using ASM I converges to the reference $\nu$ much fast than the ones of ASM II.}
    \label{Burgers_inv_nu}
\end{figure}

We point out here that in the implementation, we set $\nu = \hat{\nu} \times 0.01/\pi$, and learn the value of $\hat{\nu}$ in practice. We show representative cases for the convergences of the learned values of $\nu$ by using ASM I and ASM II in Fig. \ref{Burgers_inv_nu}. We observe from which that it is more efficient and stable by using ASM II compared with using ASM I. 
Furthermore, we run the code 10 times and present the average number of added points and the relative error for the value of $\nu$ in Table \ref{tab:burgers_inv_nu} showing that it is more accurate by using ASM II. We conclude in this subsection that the adaptive sampling method based on both residual and the gradient of the solution improves the stability and efficiency of the training procedure and the accuracy of the prediction.

\begin{table}[htbp]
\centering
 \begin{tabular}{c|cc}
  \hline \hline
  ~ & ARM I & ARM II  \\ 
  \hline
 Average \# of Added points    &  36 & 31   \\
 Relative error for $\nu$    &  $(8.95 \pm 4.95) \%$ & $(5.60 \pm 2.20) \%$  \\
  \hline \hline
 \end{tabular}
  \caption{\label{tab:burgers_inv_nu}
  Results for the inverse problem of Burgers equation. For ASM I, $m_{R} = 2$ while for ASM II, $m_{R} = 1, \; m_{\nabla u} =1$: Average numbers of added points and relative errors for the coefficient $\nu$ using ASM I and ASM II from 10 runs.}
\end{table}

\subsection{Euler equation}\label{sec:euler}
To demonstrate the present ARM for discontinuous problems, we now consider a multi-physics problem, i.e., the following one-dimensional Euler equation
\begin{align}
    \partial_t \bm{U} + \partial_x \bm{f} = 0, ~x \in [0, 1], ~ t \in [0, 2],
\end{align}
where 
\begin{equation*}
    \bm{U} =\left\{
\begin{aligned}
&\rho \\
& \rho u \\
&E
\end{aligned}
\right.,  ~\mbox{and} ~
    \bm{f} =\left\{
\begin{aligned}
&\rho u \\
&p + \rho u \\
& u (p +   E)
\end{aligned}
\right., ~
\end{equation*}
where $\rho$ is the fluid density, $u$ is the velocity, $p$ is the pressure, and $E$ is the total energy. In the present study, we consider the ideal poly-tropic gas as 
\begin{equation*}
    p = \left(\gamma - 1 \right) \left( E - \frac{1}{2} \rho u^2\right),
\end{equation*}
where $\gamma = 1.4$ is the adiabatic index. We consider the Riemann problem with a initial shock at $x = 0.5$, i.e., the initial condition is given by 
\begin{align}
    \left(\rho_L, u_L, p_L \right) = (1.4, 0.1, 1.0), ~ \left(\rho_R, u_R, p_R \right) = (1.0, 0.1, 1.0),
\end{align}
where $L$ and $R$ denote the left and right parts divided by the shock, respectively. In addition, the Dirichlet boundary conditions are imposed on the left and right boundaries. The exact solutions for this problem read as follows:
\begin{equation*}
    \rho(x, t) = \left\{ 
    \begin{aligned}
    &1.4, ~ x < 0.5 + 0.1t, \\
    & 1.0, ~ x > 0.5+0.1t,
    \end{aligned}
    \right. ~ u(x, t) = 0.1, ~ p(x, t) = 1.0,
\end{equation*}
    
\begin{figure}[htbp]
    \centering
    \subfigure[]{\label{eulera}
    \includegraphics[width=0.45\textwidth]{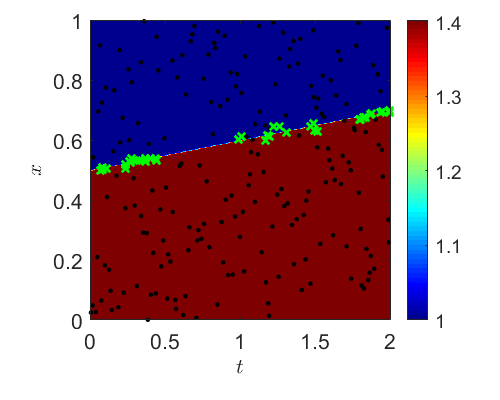}}
    % \subfigure[]{\label{eulerc}
    % \includegraphics[width=0.45\textwidth]{./Figs/euler_err.png}}
    \subfigure[]{\label{euler_cmpa}
    \includegraphics[width=0.45\textwidth]{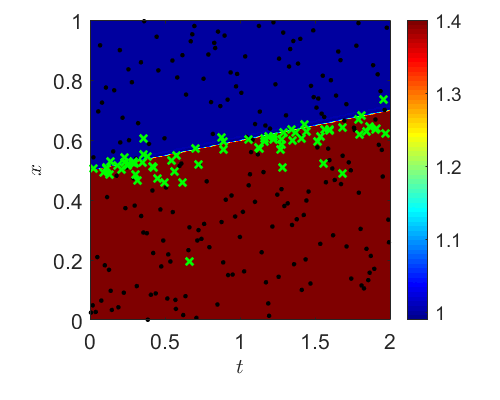}}
    \subfigure[]{\label{euler_cmpb}
    \includegraphics[width=0.45\textwidth]{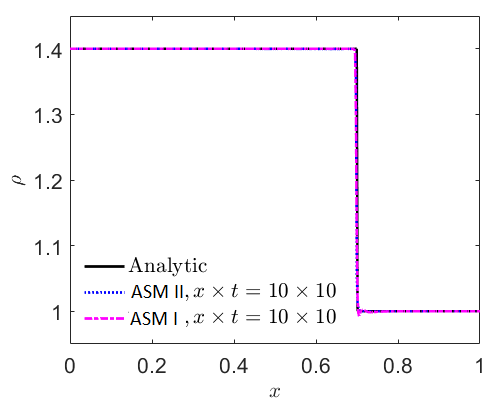}}
    \subfigure[]{\label{eulerb}
    \includegraphics[width=0.45\textwidth]{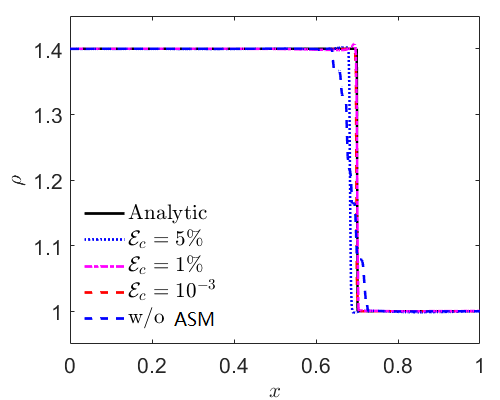}}
    \caption{
    Results for the Euler equation, $m_{R} = 1$ with ASM I while $m_{R} = m_{\nabla u} = 1$ with ASM II. 
    (a) Distribution of the residual points for solving the Euler equation using ASM II (The total number of the residual points is 236). Black dots: initial residual points, Green cross: added residual points. 
    (b) Distribution of the residual points using ASM I (The total number of the residual points is 271). 
    (c) Comparison between the solutions for $\rho$ obtained by using ASM I abd ASM II at $t = 2$.
    (d) Comparison between the solutions for $\rho$ obtained by using ASM II and the analytic solution at $t = 2$. Blue dashed line: Average value for 10 runs with 236 random residual points without the ASM. 
    % (c) History of the relative $L_2$ for $\rho$.Black triangle: mean value for 10 runs without the ARM based on 236 random residual points; Black solid line: error bar for 10 runs.
    }\label{euler}
\end{figure}

We first employ 200 random initial residual points to train the PINNs, which have 6 layers with 20 neurons per layer. In addition, the learning rate for the Adam is set as $10^{-3}$. The whole domain is divided into $x \times t = 10 \times 10$   subdomains, and we use $x \times t = 10 \times 10$ Gaussian quadrature points in each subdomain. We use ASM II and set $\mathcal{M} = 40$, $m_{R}=1$, $m_{\nabla u} = 1$, $\mathcal{E}_c = 10^{-3}$. The total number of added residual points is 236. As we can see in Fig. \ref{eulera}, all the added points are near the discontinuity, which means the ASM II cannot only roughly detect the location of the discontinuity but also resolve the issue of the discontinuity automatically, resulting in an accurate prediction for the density $\rho$ (see Fig. \ref{eulerb}). 
To further demonstrate the effectiveness of the gradient of the solution for detecting the discontinuity, here we conduct a comparison case in which we do not take the advantage of the gradient, i.e., we use ASM I with parameters $\mathcal{M} = 40$, $m_{R}=1$, and $\mathcal{E}_c = 10^{-3}$. The numbers of subdomains as well as the quadrature points are kept the same as the previous case. The number of added points is 71 (Fig. \ref{euler_cmpa}), which is about twice of the case with ASM II. We observe that most of the added points are not near the shock, leading to the slow convergence in this case. 
In addition, when we compare the predicted density obtained by using ASM I with the analytic solution, small fluctuations can be observed on the right side of the shock, which does not exist in the case when using ASM II (see Fig. \ref{euler_cmpb}). 
Furthermore, we also show the results for different $\mathcal{E}_c$ (i.e., $5\%$ and $1\%$ in Fig. \ref{eulerb}). We see that the location of the predicted discontinuity for $\mathcal{E}_c = 5\%$ differs from the analytic result. As we decrease $\mathcal{E}_c$ to $1\%$, the discontinuity is well captured although small oscillations are observed near the left part of the discontinuity. If we further use a smaller $\mathcal{E}_c = 10^{-3}$, the discontinuity becomes perfectly captured. Moreover, the oscillations near the discontinuity are also eliminated.

%For the purpose of comparison, we plot the results without the ARM using ten different groups of datasets. Each group has 236 random residual points. The averaged density from 10 runs at $t = 2$ is displayed in Fig. \ref{eulerb} indicating that the issue of the discontinuity is not resolved at all due to the lack of residual points.  The above results demonstrate the superiority of the proposed ARM. We also test the effect of the number of subdomains and quadrature points, and the conclusion is the same as the previous cases.

\begin{figure}[htbp]
    \centering
    % \subfigure[]{\label{euler_cmpa}
    % \includegraphics[width=0.45\textwidth]{./Figs/euler_data_nog.png}}
    % \subfigure[]{\label{euler_cmpb}
    % \includegraphics[width=0.45\textwidth]{./Figs/euler_cmp_nog.png}}
    \subfigure[]{\label{euler_cmpc}
    \includegraphics[width=0.45\textwidth]{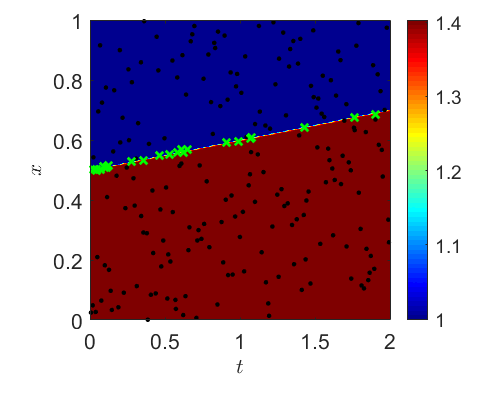}}
    \subfigure[]{\label{euler_cmpd}
    \includegraphics[width=0.45\textwidth]{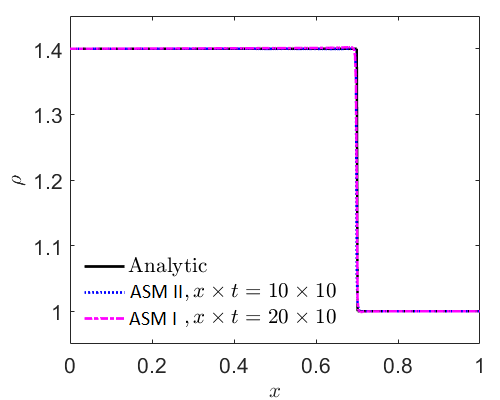}}
    % \subfigure[]{\label{eulerc}
    % \includegraphics[width=0.45\textwidth]{./Figs/euler_err.png}}
    \caption{\label{euler_cmp}
    Results for the Euler equation using sub-domains $x\times t = 20\times 10.$ $m_{R} = 1$ with ASM I while $m_{R} = m_{\nabla u} = 1$ with ASM II. 
    (a) Distribution of the residual points for the Euler equation using the PINN. Black dots: initial residual points, Green cross: added residual points. 
    (b) Comparison between the solutions for $\rho$ obtained by using ASM and the analytic solution at $t = 2$.
    % (c) History of the relative $L_2$ for $\rho$.Black triangle: mean value for 10 runs without the ARM based on 236 random residual points; Black solid line: error bar for 10 runs.
    }
\end{figure}

We also investigate how the number of sub-domains affects the performance of ASM I by increasing the number of subdomains to $x \times t = 20 \times 10$ while fixing the number of quadrature points. Almost all the added points (total number is 223) are very close to the shock, which helps resolve the sharp interface (Fig. \ref{euler_cmpc}) issue. In addition, the predicted density profile is well in agreement with the analytic solution (Fig. \ref{euler_cmpd}), which is reasonable because the computational error for the mean residual becomes smaller in smaller subdomains with the same quadrature points, which aids to detect the discontinuity. The above results indicate that the gradient of solution is capable of detecting the discontinuity efficiently.

%\subsection{Two-dimensional Poisson equation}

\subsection{Poisson equation over an L-shape domain}\label{sec:Poisson:Lshape}
Consider the following Poisson equation over an L-shaped domain $\Omega = [-1,1]^2\setminus [0,1]^2$:
\begin{align*}
    -\Delta u(x,y) &= 1,\quad (x,y)\in \Omega, \\
    u(x,y) &= 0, \quad  (x,y)\in \partial \Omega.
\end{align*}
We initially sample 400 randomly distributed residual points for all cases. For the BCs, we sample 120 randomly distributed points for either ASM I or ASM II.
The architecture of the NNs is again $4\times 40$. Similar as the case of the forward problem of Burgers equation, at each iteration, we first train the neural network using Adam optimizer with max epoch 10000 and learning rate $10^{-3}$, then apply the L-BFGS optimizer if the loss is less than $10^{-3}$. For the division of the sub-domains, we initially divide a square domain $[-1,1]^2$ into $x \times y = 25 \times 25$ uniform subdomains, and then delete the sub-domains that are in the domain $[0,1]^2$. The number of the Gaussian quadrature points in each domain is set as $12 \times 12$. The other parameters are set to be $\mathcal{M} = 80$, $\mathcal{E}_c = 3\%$, and $m_{R} = 2$ for ASM I while $m_{R} = 1, \; m_{\nabla u} =1$ for ASM II. 

\begin{figure}[htbp]
    \centering
    \subfigure[Reference solution]{\label{u_hp}
    \includegraphics[width=0.45\textwidth]{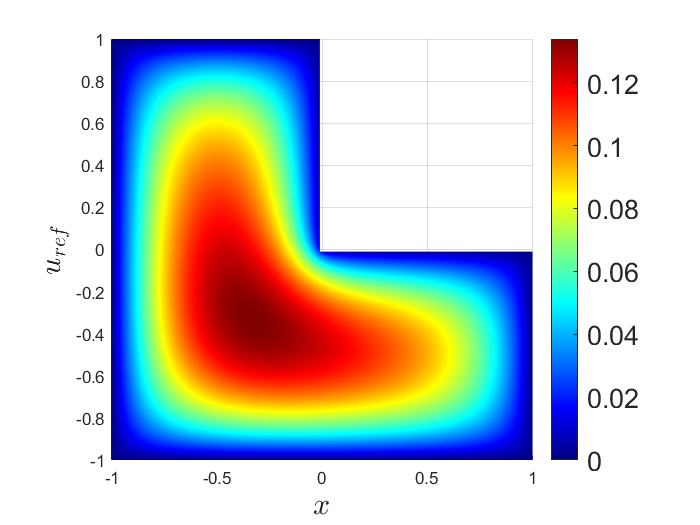}}
    \subfigure[ASM I: $E\{|u_{ref} - u_{NN}|\}$]{\label{err_res_avg}
    \includegraphics[width=0.45\textwidth]{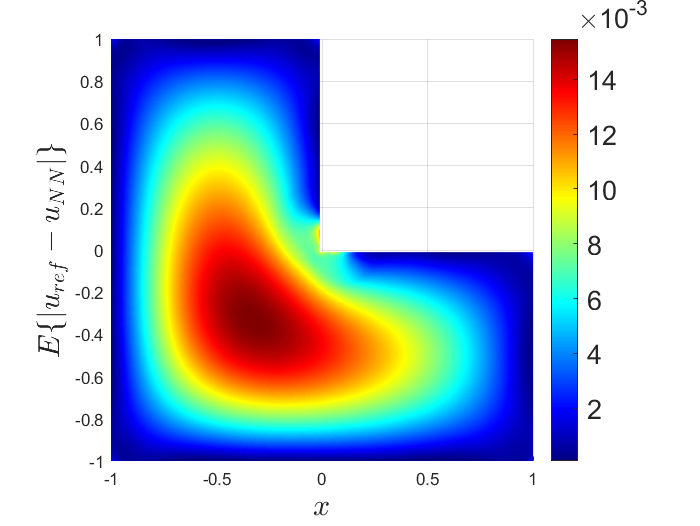}}
    \subfigure[ASM II: $E\{|u_{ref} - u_{NN}|\}$]{\label{err_res_g_avg}
    \includegraphics[width=0.45\textwidth]{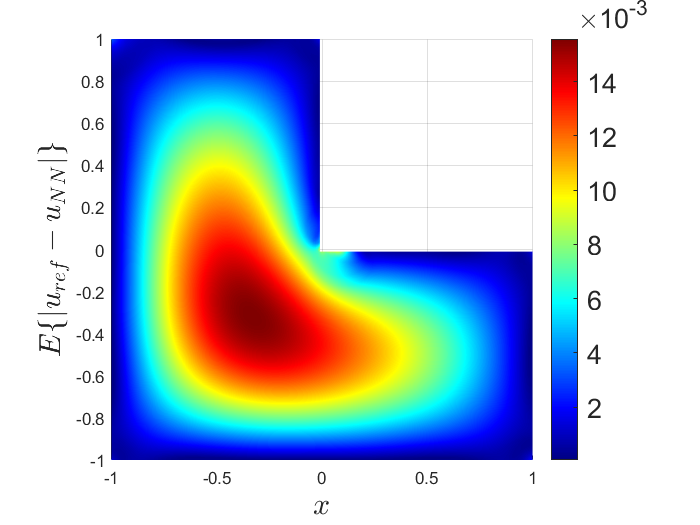}}
    \subfigure[ASM III: $E\{|u_{ref} - u_{NN}|\}$]{\label{err_res_g_bc_avg}
    \includegraphics[width=0.45\textwidth]{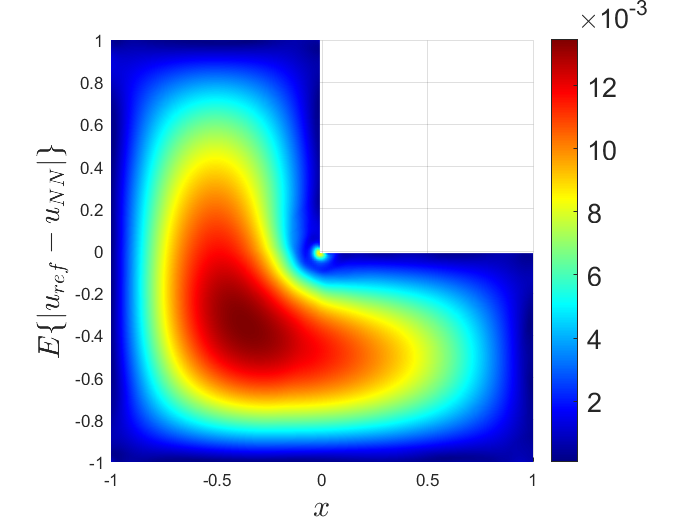}}
    \caption{Results for the Poisson equation over the L-shape domain. $m_{R} = 2$ for ASM I, $m_{R} = m_{\nabla u} = 1$ for ASM II, $m_{R} = m_{\nabla u} = 1, m_{BC} = 2$ for ASM III. (a) Reference solution obtained by using spectral element method. (b) Mean absolute error from 10 runs using ASM I. (c) Mean absolute error from 10 runs using ASM II. (d) Mean absolute error from 10 runs using ASM III, i.e., the adaptive sample technique is also applied for the boundary conditions.
    }\label{Possion_Lshape_err_avg}
\end{figure}

To investigate the stability and predictive accuracy of the proposed ASM I and II, we implement the simulation 10 times for each method. 
We show the reference solution as well as the mean values of the point-wise absolute errors for ASM I and ASM II in Figs. \ref{u_hp}-\ref{err_res_g_avg}. There is not too much difference between the results obtained by using ASM I and ASM II.
We also show the results including the average number of added residual points, Relative $L^2$ and $L^\infty$ errors by using ASM I and ASM II in Table \ref{tab:L-shape} (see the second and third rows). We observe that the number of added residual points and the relative $L^2$ errors for the cases of ASM I and ASM II are very close to each other, only a weak advantage of ASM II is observed in the sense of accuracy and stability. 
The reason is that both the mean value of the point-wise error and the $L^2$ error are in the average sense. 
However, we observe from the $L^\infty$ error that the ASM II algorithm delivers much better predictive accuracy and stability compared with the ASM I algorithm. This means that the gradient of the solution helps a lot in the adaptive sampling procedure.

\begin{figure}[htbp]
    \centering
    \subfigure[ASM I: $|u_{ref} - u_{NN}|$]{\label{err_res}
    \includegraphics[width=0.45\textwidth]{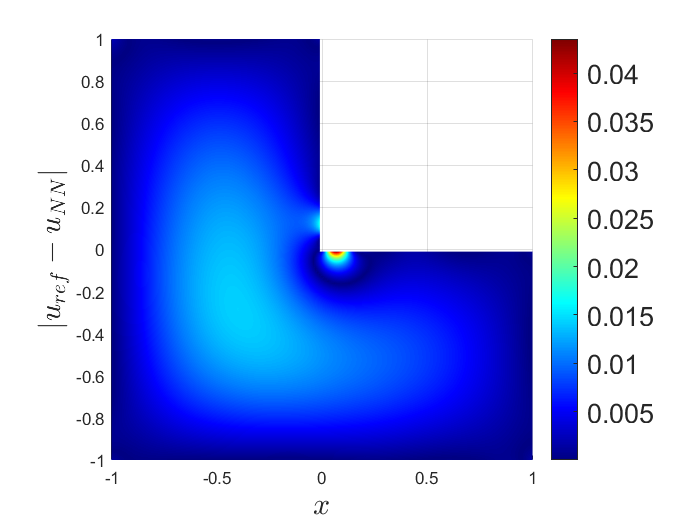}}
    \subfigure[ASM II: $|u_{ref} - u_{NN}|$]{\label{err_res_g}
    \includegraphics[width=0.45\textwidth]{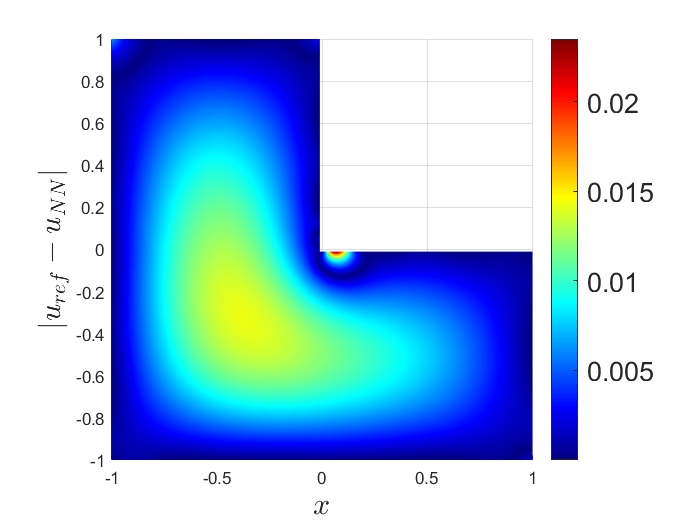}}
    \caption{Results for the Poisson equation defined on the L-shape domain. (a) A representative absolute error by using ASM I. (b) A representative absolute error by using ASM II.
    }\label{Possion_Lshape_err}
\end{figure}

For the L-shape Poisson problem, we see that the sharpness is near the corner. Moreover, unlike the Burgers equation whose solution sharpness is far away from boundaries, in the case of the L-shape Poisson problem, the sharpness of the solution is near boundaries. We found that this is one of the main reasons resulting the instability of the proposed methods. Two representative cases for ASM I and ASM II are shown in Fig. \ref{Possion_Lshape_err}. Observe that large errors are presented near boundaries. 

Therefore,  we further employ the adaptive sampling method for data points for the boundary conditions. The detail of the algorithm are given in Algorithm \ref{alg:ASMIII}. 
We point out here that for the sake of fair comparison, we initially sample 80 points for BCs and set $m_{BC} = 2$. Also, we set a minimum number of iterations to be 9 in order to sufficiently train the network.
The result of the mean value of the point-wise absolute error is shown in Fig. \ref{err_res_g_bc_avg} while the relative $L^2$ and $L^\infty$ errors are shown in the forth row in Table \ref{tab:L-shape}, showing that the ASM III algorithm has the supreme accuracy and stability.

% \begin{figure}
%     \centering
%     \subfigure[]{\label{u_hp}
%     \includegraphics[width=0.45\textwidth]{./Figs/u_ref.png}}
%     \subfigure[]{\label{u_res}
%     \includegraphics[width=0.45\textwidth]{./Figs/u_res.png}}
%     \subfigure[]{\label{u_res_g}
%     \includegraphics[width=0.45\textwidth]{./Figs/u_res_g.png}}
%     \subfigure[]{\label{u_res_g_bc}
%     \includegraphics[width=0.45\textwidth]{./Figs/u_res_g_bc.png}}
%     \caption{\label{Possion_Lshape_u}
%     }
% \end{figure}

\begin{table}[htbp]
\centering
 \begin{tabular}{c|ccc}
  \hline \hline
  ~ & ASM I & ASM II & ASM III (for BCs) \\ 
  \hline
 Average \# of Added points    &  57 & 60 & 20  \\
 Relative $L^2$ error  &  $(10.57 \pm 1.92) \%$ & $(10.39 \pm 1.12) \%$  & $(8.60 \pm 0.91) \%$\\
 Relative $L^\infty$ error  &  $(25.38 \pm 14.44) \%$ & $(14.48 \pm 3.08) \%$  & $(9.80 \pm 1.01) \%$\\
  \hline \hline
 \end{tabular}
  \caption{\label{tab:L-shape}
  Results for the Poisson problem over the L-shape domain. $m_{R} = 2$ for ASM I, $m_{R} = m_{\nabla u} = 1$ for ASM II, $m_{R} = m_{\nabla u} = 1, m_{BC} = 2$ for ASM III: Average numbers of added points and relative $L^2$ and $L^\infty$ errors using ASM I, ASM II and ASM III. We run the code 10 times for each case.}
\end{table}

\begin{algorithm} [H]
%    \SetAlgoNoLine  %去掉之前的竖线 
\caption{ASM III: Residual/gradient/BC-based adaptive sampling method}\label{alg:ASMIII}
\begin{enumerate}
\item[Step 1] Train the PINNs based on the residual point set $\mathcal{T}$ and the boundary point set $\mathcal{B}$;
\item[Step 2] Divide the computational domain into $\bm{S}$  subdomains $\Omega_i, i =1,\ldots, {\bm S}$, and compute 
$$\mathcal{E}_i = \int_{\Omega_i} |R| d{\bm x}, \quad
\mathbf{E}_i = \int_{\Omega_i} |\nabla u| d{\bm x}
$$
using the Gaussian quadrature method for each sub-domain;
using the Gaussian quadrature method for each sub-domain;
\item[Step 3] Let $\mathcal{E}_{max} = \max_{i = 1, \ldots {\bm S}}\{\mathcal{E}_i\}$ and $\mathbf{E}_{max} = \max_{i = 1, \ldots {\bm S}}\{\mathbf{E}_i\}$. If $\mathcal{E}_{max} < \mathcal{E}_c$, iteration stops. Otherwise, randomly sample $\mathcal{M}$ data in the sub-domain which has the largest $\mathcal{E}$ and $\mathcal{M}$ data in the sub-domain which has the largest $\mathbf{E}$, and then compute the corresponding residuals and gradients for selected points, respectively. Add $m_{R}$ points which have the largest $|R|$ and $m_{\nabla u}$ points which have the largest $|\nabla u|$ in the residual point set .
\item[Step 4] Randomly sample $\mathcal{M}_{BC}$ boundary points and compute the absolute value of $u$ at these boundary points. Then add $m_{BC}$ points which have the largest $|u|$ in the boundary point set $\mathcal{B}$.
\item[Step 5] Repeat Step 1 - 4.
\end{enumerate}
\end{algorithm}

\subsection{High-dimensional Poisson equation}\label{sec:Poisson:high}
Consider the following high-dimensional Poisson equation:
\begin{align*}
    -\Delta u({\bm x}) &= f({\bm x}),\quad {\bm x}\in \Omega:=[-1,1]^d, \\
    u({\bm x}) &= 0, \quad  {\bm x}\in \partial \Omega.
\end{align*}
Same as that used in \cite{gao2022failure}, we let the exact solution be $u({\bm x}) = e^{-10\|{\bm x}\|^2}$ and set $d = 9$. 20000 residual points and 7200 BC points are initially sampled, and then we add the residual points adaptively by using ASM I or ASM II.
The architecture of the NNs is $6\times 40$. Again, at each iteration, we first train the neural network using the Adam optimizer with max epoch 10000 and learning rate $8\times 10^{-4}$, then apply the L-BFGS optimizer if the loss is less than $10^{-4}$. 
In this case, we don't use the domain decomposition, therefore, the ASM I is equivalent to the residual adaptive refinement method proposed in \cite{lu2021deepxde}. For both ASM I and ASM II, at each iteration, we randomly sample 300000 points and determine the new added residual points accordingly. Also, we set $m_{R} = 500$ for ASM I while $m_{R} = 300, \; m_{\nabla u} = 200$ for ASM II. 
In this case, instead of training the network to a given $\mathcal{E}_c$, we train the network with a maximum iteration number 12.
%$\mathcal{E}_c = 30\%$, and 

\begin{figure}
    \centering
    \includegraphics[width = 0.6\textwidth]{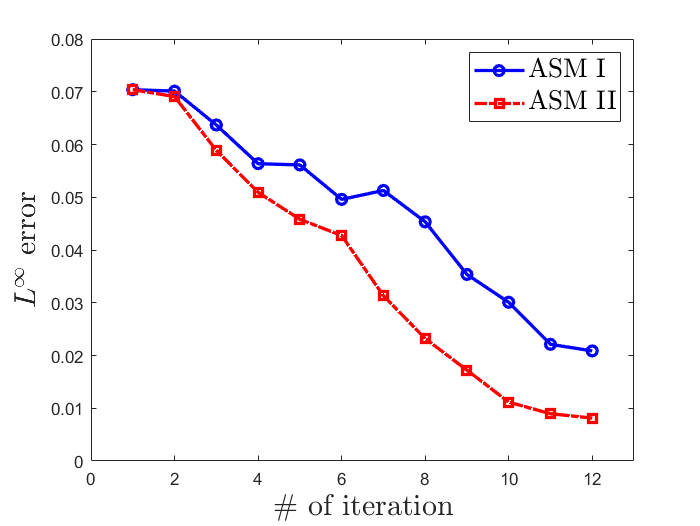}
    \caption{Convergence of the $L^\infty$ error for the high dimensional Poisson problem using ASM I and ASM II. $m_{R} = 500$ for ASM I while $m_{R} = 300, \; m_{\nabla u} = 200$ for ASM II.}
    \label{Possion_HighD_Linf}
\end{figure}

The comparison of the $L^\infty$ error between ASM I and ASM II is shown in Fig. \ref{Possion_HighD_Linf}, showing again that the error obtained by using ASM II converges fast than that obtained by using ASM I.

\section{Concluding remarks}
\label{summary}

Adaptive sampling methods (ASMs) for solving the forward and inverse PDE problems, which have sharp solutions, using the physics-informed neural networks are presented in the present study. The key point of the ASM is to add more residual points based the PDE residual and/or the gradient of the predicted solution. To be more specific: we proposed two types of ASM algorithms, the first one is residual only-based denoted by ASM I, while the other one is based on both the residual and the gradient of the solution denoted by ASM II.
The ASMs are then applied to solving several PDEs: 
\begin{enumerate}
    \item The forward and inverse problems of Burgers equation;
    \item A Riemann problem of the compressible Euler equation;
    \item Poisson equations defined on an L-shape domain and a high-dimensional hypercube domain.
\end{enumerate}
In all test cases, the ASM delivers more accurate results than PINNs with the same number of residual points. Moreover, by utilizing the gradient of the  predicted solution, the ASM is more stable and capable of detecting and refining the discontinuities efficiently as well as delivers better accuracy compared with the residual only-based ASM. Finally, we also extended the adaptive refinement strategy to the boundary conditions for the problem whose solution has sharpness near the boundary giving a more stable training process and more accurate solutions.

\section*{Acknowledgements}
Z. Mao would like to acknowledge the support of National Natural Science Foundation of China under Grant No. 12171404 and National Key R\&D Program of China under Grant No. 2022YFA1004504.
%Fundamental Research Funds for Central Universities of China under Grant No. 20720210037.
X. Meng would like to acknowledge the support of the National Natural Science Foundation of China (No. 12201229)

%Set $\varepsilon=1/10$.

%% The Appendices part is started with the command \appendix;
%% appendix sections are then done as normal sections
%% \appendix

%% \section{}
%% \label{}

%% References
%%
%% Following citation commands can be used in the body text:
%% Usage of \cite is as follows:
%%   \cite{key}          ==>>  [#]
%%   \cite[chap. 2]{key} ==>>  [#, chap. 2]
%%   \citet{key}         ==>>  Author [#]

%% References with bibTeX database:

\bibliographystyle{elsarticle-num}
\biboptions{sort&compress}
\bibliography{adaptive.bib}

%% Authors are advised to submit their bibtex database files. They are
%% requested to list a bibtex style file in the manuscript if they do
%% not want to use model1-num-names.bst.

%% References without bibTeX database:

% \begin{thebibliography}{00}

%% \bibitem must have the following form:
%%   \bibitem{key}...
%%

% \bibitem{}

% \end{thebibliography}

\end{document}